\documentclass[a4paper,11pt]{amsart}
\usepackage{amssymb}
\usepackage{amscd}

\newtheorem{theorem}{Theorem}[section]

\newtheorem{proposition}[theorem]{Proposition}
\newtheorem{corollary}[theorem]{Corollary}

\theoremstyle{definition}

\theoremstyle{remark}

\usepackage{graphicx}
\usepackage{amsmath}
\usepackage{amsfonts}
\usepackage{amssymb}
\newcommand{\be}{\begin{equation}}

\newcommand{\ee}{\end{equation}}

\newcommand{\Om}{\Omega}

\newcommand{\K}{\Upsilon}

\newcommand{\om}{\omega}

\newcommand{\cB}{{\mathcal B}}

\newcommand{\si}{\sigma}

\newcommand{\ba}{\begin{array}}

\newcommand{\ea}{\end{array}}

\newcommand{\beq}{\begin{eqnarray}}

\newcommand{\eeq}{\end{eqnarray}}

\newtheorem{lm}{lemma}

\newtheorem{thee}{theorem}

\newtheorem{proo}{proposition}

\newtheorem{co}{corollary}

\newtheorem{rem}{remark}

\newtheorem{deff}{definition}

\newcommand{\bd}{\begin{deff}}

\newcommand{\ed}{\end{deff}}

\newcommand{\bl}{\begin{lm}}

\newcommand{\el}{\end{lm}}

\newcommand{\bp}{\begin{proo}}

\newcommand{\ep}{\end{proo}}

\newcommand{\bt}{\begin{thee}}

\newcommand{\et}{\end{thee}}

\newcommand{\bc}{\begin{co}}

\newcommand{\ec}{\end{co}}

\newcommand{\brm}{\begin{rem}}

\newcommand{\erm}{\end{rem}}

\usepackage{t1enc}

\def\Bbb{\mathbb}
\def\Cal{\mathcal}

\newcommand{\newc}{\newcommand}

\let\ccdot\cdot
\def\cdot{\hbox to 2.5pt{\hss$\ccdot$\hss}}

\newc{\aR}{\mbox{\boldmath{$ R$}}}
\newc{\aS}{\mbox{\boldmath{$ S$}}}
\newc{\aT}{\mbox{\boldmath{$ T$}}}
\newc{\aW}{\mbox{\boldmath{$ W$}}}

\newc{\aK}{\mbox{\boldmath{$ K$}}}
\newc{\aL}{\mbox{\boldmath{$ L$}}}

\newcommand{\ce}{{\Cal E}}

\newcommand{\cq}{{\Cal Q}}

\newcommand{\ct}{{\Cal T}}

\usepackage{amssymb}
\usepackage{amscd}

\newcommand{\nd}{\nabla}

\newcommand{\Rho}{{\mbox{\sf P}}}
\newcommand{\Up}{\Upsilon}

\let\hash=\sharp  

\newcommand{\fb}  {\mbox{$                                      
\begin{picture}(9,8)(1.6,0.15)
\put(1,0.2){\mbox{$ \Box \hspace{-7.8pt} /$}}
\end{picture}$}} 

\newcommand{\wh}{\widehat}

\newcommand{\cV}{{\Cal V}}

\newcommand{\Pa}{{\Bbb I}}

\newcommand{\nn}[1]{(\ref{#1})}

\newcommand{\bg}{\mbox{\boldmath{$ g$}}}

\let\t=\tau

\newcommand{\V}{{\mbox{\sf P}}}                   
\newcommand{\J}{{\mbox{\sf J}}}
\newc{\obstrn}[2]{B^{#1}_{#2}}

\newcommand{\rpl}                         
{\mbox{$
\begin{picture}(12.7,8)(-.5,-1)
\put(0,0.2){$+$}
\put(4.2,2.8){\oval(8,8)[r]}
\end{picture}$}}

\newcommand{\lpl}                         
{\mbox{$
\begin{picture}(12.7,8)(-.5,-1)
\put(2,0.2){$+$}
\put(6.2,2.8){\oval(8,8)[l]}
\end{picture}$}}

\usepackage{ifthen}

\newc{\tensor}[1]{#1}
\newc{\Mvariable}[1]{\mbox{#1}}
\newc{\down}[1]{{}_{#1}}
\newc{\up}[1]{{}^{#1}}

\newc{\JulyStrut}{\rule{0mm}{6mm}}
\newc{\midtenPan}{\mbox{\sf S}}
\newc{\midten}{\mbox{\sf T}}
\newc{\midtenEi}{\mbox{\sf U}}
\newc{\ATen}{\mbox{\sf E}}
\newc{\BTen}{\mbox{\sf F}}
\newc{\CTen}{\mbox{\sf G}}

\def\sideremark#1{\ifvmode\leavevmode\fi\vadjust{\vbox to0pt{\vss
 \hbox to 0pt{\hskip\hsize\hskip1em
 \vbox{\hsize3cm\tiny\raggedright\pretolerance10000
 \noindent #1\hfill}\hss}\vbox to8pt{\vfil}\vss}}}%
                        
\numberwithin{equation}{section}

\begin{document}
\renewcommand{\today}{}
\title{Almost conformally Einstein manifolds and obstructions}
\author{A. Rod Gover}

\address{Department of Mathematics\\
  The University of Auckland\\
  Private Bag 92019\\
  Auckland 1\\
  New Zealand} \email{gover@math.auckland.ac.nz}

\vspace{10pt}

\renewcommand{\arraystretch}{1}
\maketitle
\renewcommand{\arraystretch}{1.5}

\pagestyle{myheadings}
\markboth{Gover}{Almost conformally Einstein manifolds and obstructions}

\begin{abstract}
  A Riemannian or pseudo-Riemannian (or conformal) structure is
  conformally Einstein if and only if there is a suitably generic
  parallel section of a certain vector bundle -- the so-called
  standard conformal tractor bundle.  We show that this
  characterisation leads to a systematic approach to constructing
  obstructions to conformally Einstein metrics.  Relaxing the
  requirement that the parallel tractor field be generic gives a
  natural generalisation of the Einstein equations.
\end{abstract}
\thanks{The author gratefully
  acknowledges support from the Royal Society of New Zealand via
  Marsden Grant no.\ 02-UOA-108, and also the New Zealand Institute of
  Mathematics and its Applications for support via a Maclaurin
  Fellowship. The work was partially done while the author was visiting
  the Institute for Mathematical Sciences, National University of
  Singapore in 2004. The visit was supported by the Institute.}


\maketitle
\section{Introduction}

In these partly expository notes we review some recent results
concerning obstructions for metrics (on manifolds of dimension $n\geq
3$) to be conformally Einstein and discuss the relationship between
these and a generalisation of the Einstein condition.  Given an
initial metric we show that the equations for a conformally related
metric to be Einstein, when prolonged and written as a first order
system, give a conformally invariant connection on a certain vector
bundle. This is the so-called (standard conformal) tractor connection
and bundle respectively.  A solution of the original equations, that
is an Einstein metric, is then equivalent to a suitably generic
parallel section of the tractor bundle. Relaxing the condition that
the parallel tractor be generic then leads to an obvious
generalisation of the Einstein equations. These generalised
structures, viz.\ Riemannian or pseudo-Riemannian manifolds equipped
with a parallel standard tractor, are termed almost Einstein
structures.  A closely related term is almost conformally Einstein: we
say a Riemannian or pseudo-Riemannian manifold is almost conformally
Einstein if on an open dense subset it is conformally related to
Einstein metric.  In sections \ref{two} and \ref{three} we include
some new (although elementary) analysis of the zero set of the scaling
function on almost Einstein spaces, that is the set where the Einstein
metric is singular.  There is no attempt to be complete in this
treatment. Rather it is intended to merely point out some of the most
obvious properties of parallel tractors.

The problem of finding necessary and sufficient conditions for a
Riemannian or pseudo-Riemannian manifold to be locally conformally
related to an Einstein metric has been studied for almost 100 years.
Early results date back to the work of Brinkman \cite{B1,B2} and
Schouten \cite{S}. Substantial progress was made by Szekeres in the
1963 \cite{Sz} and then Kozameh, Newman and Tod (KNT) \cite{CTP} in
1980's.  One approach is to seek invariants, polynomial in the
Riemannian curvature and its covariant derivatives, that give a sharp
obstruction to conformally Einstein metrics in the sense that they
vanish if and only if the metric concerned is conformally related to
an Einstein metric. For example in dimension 3 it is well known that
this problem is solved by the Cotton tensor, which is a certain tensor
part of the first covariant derivative of the Ricci tensor. This
tensor is also a sharp obstruction to local conformal flatness. So
locally 3-manifolds are conformally Einstein if and only if they are
conformally flat.  In \cite{CTP} KNT described conformal invariants
that gave sharp obstructions on 4 manifolds given a restriction that
the class of metrics to be considered is suitably generic. One
component of the KNT system is the Bach tensor. For manifolds of
higher even dimension there is a natural analogue of the Bach tensor;
this is the Fefferman-Graham obstruction tensor. This trace-free,
divergence-free (density valued) 2-tensor is of considerable current
interest, in part due to its relationship to Branson's Q-curvature
\cite{GrH}. In the original work of Fefferman and Graham this tensor
$\cB$ arose as an obstruction to their ambient metric construction
\cite{FeffGr}.  They observed that for conformally Einstein manifolds
the ambient construction works to all orders and hence $\cB$ vanishes
identically (see \cite{GrH} where the equivalent argument in terms of
the Poincar\'e metric is given explicitly).

One point of this article is to illustrate that the tractor bundle and
the relationship of parallel tractors to Einstein metrics leads to a a
uniform treatment of these obstructions. We review here aspects of the
recent work \cite{GNur} of the author with Nurowski where it is
observed that many of the classical obstructions to conformally
Einstein metrics, including the KNT invariants, can be recovered via
the integrability conditions for a parallel tractor. In fact the
treatment via tractors leads to new obstructions and in particular a
system of invariants that gives a sharp obstruction to conformally
Einstein metrics for the class of metrics that are weakly generic
(this means that, viewed as a bundle map $TM\to\otimes^3 TM$, the Weyl
curvature is injective). This is a broader class of metrics than
treated previously and in particular extends considerably the recent
results of Listing \cite{Listing}. Following \cite{GoPetob} we then
observe that essentially the same ideas lead to a simple and direct
proof that the Fefferman-Graham tensor $\cB$ vanishes for conformally
Einstein metrics. Aside from these new results for Riemannian and
pseudo-Riemannian geometry a point that should be made is that the
ideas here for the construction and study of obstructions via parallel
tractor fields adapt easily to other equations and to other structures. An
obvious example is projective structures where one could follow these
ideas to give obstructions to Ricci flat scales. There is also
evidence that the related ideas will have a useful role in CR geometry
\cite{CapGoFeff}, \cite{GLeit}.

Of course, by continuity, conformal invariants which vanish for
conformally Einstein metrics also vanish for almost conformally
Einstein metrics. Thus in any case where such invariants give a sharp
obstruction the equations they determine are naturally viewed as the
equations for almost conformally Einstein metrics.

\section{Einstein metrics and the tractor bundle} \label{two}

Let $M$ be a smooth manifold, of dimension $n\geq 3$, equipped with a
Riemannian or pseudo-Riemannian metric $g_{ab}$. Here and throughout
we employ Penrose's abstract index notation \cite{ot} and indices
should be assumed abstract unless otherwise indicated. We write
$\ce^a$ to denote the space of smooth sections of the tangent bundle
on $M$, and $\ce_a$ for the space of smooth sections of the cotangent
bundle.  (In fact we will often use the same symbols for the
corresponding bundles, and also in other situations we will often use
the same symbol for a given bundle and its space of smooth sections,
since the meaning will be clear by context.) We write $\ce$ for the
space of smooth functions and all tensors considered will be assumed
smooth without further comment.  An index which appears twice, once
raised and once lowered, indicates a contraction.  The metric $g_{ab}$
and its inverse $g^{ab}$ enable the identification of $\ce^a$ and
$\ce_a$ and we indicate this by raising and lowering indices in the
usual way.

Recall that the Levi-Civita connection $\nabla_a$  is the unique torsion free
connection preserving the metric $g_{ab}$. The Riemann curvature tensor
$R_{ab}{}^{c}{}_d$ is given by
$$(\nabla_a\nabla_b-\nabla_b\nabla_a)V^c=R_{ab}{}^{c}{}_d V^d 
\quad\text{
where} \quad \ V^c\in \ce^c.$$
This can be decomposed into the totally trace-free {\em Weyl curvature}
$C_{abcd}$ and the symmetric {\em
Schouten tensor} $\Rho_{ab}$ according to
$$
R_{abcd}=C_{abcd}+2g_{c[a}\Rho_{b]d}+2g_{d[b}\Rho_{a]c}.
$$
Thus $\Rho_{ab}$ is a trace modification of the Ricci tensor 
${\rm Ric}_{ab}=R_{ca}{}^c{}_b$:
$$
R_{ab}=(n-2)\Rho_{ab}+\J g_{ab}, \quad \quad \J:=\Rho^a_{~a}.
$$ 
Recall that a metric is {\em Einstein} if the Ricci tensor is pure
trace. Equivalently this means the Schouten tensor is pure trace,
$$
\Rho_{ab}-\frac{1}{n} \J g_{ab}=0.
$$

Given a (pseudo-)Riemannian metric $g$ we are interested in the
question of whether $g$ is {\em conformally Einstein}. That is whether
there is some metric $\wh{g}$, conformally related to $g$
$$
\wh{g}=e^{2\om}g   \quad \om\in \ce 
$$ which is Einstein.  
We write
$\wh{\Rho}_{ab}-\frac{1}{n} \wh{\J} \wh{g}_{ab}=0 $ where the hatted
quantities refer to the metric $\wh{g}$.  From the formula for the
Levi-Civita connection in terms of the metric we can calculate the
conformal variation of the Schouten tensor. From this we obtain that,
in terms of the metric $g_{ab}$, the condition for $ \wh{g}_{ab}$ to
be Einstein is
\begin{equation}\label{confein}
\Rho_{ab}-\nabla_a\K_b+\K_a\K_b-\frac{1}{n}Tg_{ab}=0, 
\end{equation}
where
$$
\Up:= d\om
$$ and $ T=\J-\nabla^a\K_a+\K^a\K_a $.  

As an equation on the function $\om$, the system \nn{confein} is
clearly overdetermined and we do not expect solutions in general. We
will see over the following pages that there is systematic way to
construct obstructions.  First note that the system \nn{confein} can
be linearised by the change of variables $\om \mapsto \si:=e^{-\om}$.
The resulting equivalent equation on the positive function $\si$ is
$$
{\rm TF}(\nd_a\nd_b +\Rho_{ab})\si =0,  
$$ 
where ``TF'' indicates that we take the trace-free part. Next we replace this equation with the equivalent first order system:
\begin{equation}\label{syst}
\nd_a\si-\mu_a=0, \quad  \mbox{and} \quad \nd_a\mu_b+\Rho_{ab}\si +g_{ab}\rho=0,
\end{equation}
where $\mu_a$ is a 1-form field and $\rho$ a function. 
Differentiating the second of these again and contracting yields an
equation on $\rho$:
$$
\nd_a\rho -\Rho_{ab}\mu^b=0.
$$
The system has closed up linearly. Thus the original equation is
equivalent to a connection and the equation for a parallel section of
this. Let $ \Pa:=(\si,\mu_b,\rho)\in \ce\oplus \ce_b\oplus\ce$ then
$$
 \nn{confein}  \Leftrightarrow \nd \Pa =0,
$$
with the qualification that $\si$ is a positive function, and where
\begin{equation}\label{trconn}
\nd_a
\left( \begin{array}{c}
\si\\\mu_b\\ \rho
\end{array} \right) : =
\left( \begin{array}{c}
 \nabla_a \si-\mu_a \\
 \nabla_a \mu_b+ g_{ab} \rho +\Rho_{ab}\si \\
 \nabla_a \rho - \Rho_{ab}\mu^b  \end{array} \right) .
\end{equation}
Note that the formula gives a sum of the trivial extension of the
Levi-Civita connection with a bundle endomorphism of $\ce\oplus
\ce_b\oplus\ce$ and so this is a connection on $\ce\oplus
\ce_b\oplus\ce$.  Through context, no confusion should arise from the
use of the symbol $\nd$ to denote this new connection as well as the
Levi-Civita connection.

We can be more precise concerning the relationship between
\nn{confein} and the prolonged system $\nd \Pa=0$. Note that if $\Pa$
is parallel then we recover the system \nn{syst}. From the first of
these and a trace of the second we have, respectively, $\mu_a =\nd_a
\si$ and $\rho=-\frac{1}{n}(\Delta + \J) \si$ (here
$\Delta=\nd^a\nd_a$). Thus we can say that there is a 1-1
correspondence between scales $\si$ such that $\wh{g}=\si^{-2}g$ is
Einstein and sections $\Pa=(\si,\mu_a,\rho)$ of the bundle $\ce\oplus
\ce_a\oplus \ce$ with $\si$ nowhere vanishing (note we do not need
$\si$ positive for this statement). The mapping from Einstein scales
to parallel sections of $\ce\oplus \ce_b\oplus\ce$ is given by
$\si\mapsto \frac{1}{n}D \si$ where $D\si : =(n \si, n \nd_a \si, -
(\Delta + \J) \si )$.

A main point that we wish to come to is that from the connection and
its curvature it is easy to construct obstructions to the equation
\nn{confein} for a metric to be conformally Einstein. Before we do
this let us observe the conformal invariance of the connection. Of
course the question of whether there is a metric, conformally related
to $g$, that is Einstein is tautologically an issue of conformal structure. 

Recall that a  {\em conformal
  structure\/} on $M$ of signature $(p,q)$  is an
equivalence class $[g]$ of pseudo--Riemannian metrics of
signature $(p,q)$ on $M$, with two metrics being equivalent if and
only if one is obtained from the other by multiplication with a
positive smooth function. Equivalently a conformal structure is a
smooth ray subbundle $\cq\subset S^2T^*M$ whose fibre over $x$
consists of the values of $g_x$ for all metrics $g$ in the conformal
class. In this picture a metric in the conformal class is a section of $\cq$. 

We can view $ \cq$ as a principal bundle
$\pi:\cq\to M$ with structure group $\Bbb R_+$, and so there are
natural line bundles on $ (M,[g])$ induced from the irreducible
representations of ${\Bbb R}_+$. For $w\in {\Bbb R}$, we write $\ce[w]$
for the line bundle induced from the representation of weight $-w/2$
on ${\Bbb R}$ (that is, ${\Bbb R}_+ \ni x\mapsto x^{-w/2}\in {\rm
End}(\Bbb R)$). Clearly the fibres inherit an ordering from $\Bbb R$.
A section of $E[w]$ corresponds to a real-valued function $ f$ on
$\Cal Q$ with the homogeneity property $f(x, \Omega^2 g)=\Omega^w
f(x,g)$, where $\Omega$ is a positive function on $M$, $x\in M$, and
$g$ is a metric from the conformal class $[g]$. 
 Given a vector bundle $\cV$ or
section space thereof we write $\cV[w]$ as a shorthand for $\cV\otimes \ce[w]$.

There is a tautological function $\bg$ on $\cq$ taking values in
$S^2T^*M$, namely the function which assigns to the point $(x,g_x)\in
\cq$ the metric $g_x$ at $x$.  This is homogeneous of degree 2, since
$\bg (x,s^2 g_x) =s^2 g_x$ and so $\bg$ is equivalent to a section of
$\ce_{(ab)}[2]$ that we denote by the same symbol and term the {\em
conformal metric}. (Note $(a\cdots b)$ means the symmetric part over
the enclosed indices.) Then a metric $g$ from the conformal class is
determined by  a non-vanishing section $\si$ of $\ce[1]$ (a so-called {\em
conformal scale}) via the equation $g=\si^{-2}\bg$. 
 The conformal metric $\bg_{ab}$ has an inverse
$\bg^{ab}$ in $\ce^{(ab)}[-2]$.

 By using density-valued tensor fields and $\bg_{ab}$ and its inverse
to raise and lower indices and so forth we can avoid carrying around
conformal factors in calculations. Otherwise the calculations are
almost formally identical to the calculations where one instead picks
a ``reference metric''. For example in these terms the conformal scale
$\si\in \ce[1]$ gives an Einstein metric if and only if it solves the 
equation 
\begin{equation}\label{siver}
{\rm TF}(\nd_a\nd_b +\Rho_{ab})\si =0. 
\end{equation}
Here the Schouten tensor $\Rho$ and the Levi-Civita connection
$\nabla$ are constructed from any metric $g$ from the conformal class
and $\nd$ is extended to act on densities (and density valued tensors)
through the trivialisation of the density bundles given by the choice
$g$.  It is easily verified explicitly that this equation is {\em
  conformally invariant}. Indeed the left-hand-side is conformally
invariant; if $\wh{g}$ is another metric from the conformal class then
we have $ {\rm TF}(\wh{\nd}_a\wh{\nd}_b +\wh{\Rho}_{ab})\si={\rm
  TF}(\nd_a\nd_b +\Rho_{ab})\si .  $ Thus if we begin with \nn{siver},
in our argument above, then the corresponding connection \nn{trconn}
is conformally invariant. This is the so-called (conformal) tractor
connection \cite{Cartan,Thomas,BEG} (which is an induced connection
equivalent to the normal conformal Cartan connection -- see
\cite{Cap-Gover}). Note that the connection really acts on the
(standard) {\em tractor bundle} $\ct $ (or $\ce^A$ in abstract index
notation) which may be defined as the quotient of $J^2\ce[1]$ by the
image of $\ce_{(ab)_0}[1]$ in $J^2\ce[1]$ through the jet exact
sequence at 2-jets. Here $\ce_{(ab)_0}[1]$ means the trace-free
symmetric part of $\ce_{ab}\otimes \ce[1]$. By construction then there
is a canonical surjection $\ct\to J^1\ce[1]$. Composing on the left
with the jet projection $J^1\ce[1]\to \ce[1]$ we obtain a canonical
map $X_A:\ce^A\to \ce[1]$.  The invariant version of the operator $D$
above (that we again denote by $D$) is just the composition involving
the universal 2-jet operator $j^2:\ce[1]\to J^2\ce[1]$ followed by the
canonical projection $J^2\ce[1]\to \ct$.  Via a choice of metric $g$,
and the Levi-Civita connection it determines, we may use the
\underline{formula} given for ($1/n$ times) $D$ above, viz $\si\mapsto
( \si, \nd_a \si, - \frac{1}{n}(\Delta + \J) \si )$ to split the
tractor bundle $[\ce^A]_g=\ce[1]\oplus \ce_a[1]\oplus \ce[-1]$. In the
discussions that follow we will often use the splitting determined by
some choice of metric without explicit comment. In terms of such a
splitting the conformally invariant tractor connection is then exactly
given by the formula \nn{trconn}. This connection preserves a
conformally invariant signature $(p+1,q+1)$ {\em tractor metric} $h$
which, in terms of the splitting mentioned, is given by
$h(V,V)=2\si\rho + \bg^{ab}\mu_a\mu_b$, if $[V]_g=(\si,\mu_a,\rho)$.
For further details on the tractor calculus from this point of view
see \cite{luminy}. The relationship between parallel tractors and
conformally Einstein metrics, while implicit in \cite{BEG}, was
probably first observed and treated in detail by Paul Gauduchon in \cite{Gau}.  Let us
summarise our observations and results as a theorem (cf.\ theorem 3.1
in \cite{GNur}).
\begin{theorem}\label{cein}
On a conformal manifold $(M,[g])$ there is a 1-1 correspondence
between conformal scales $\si\in \ce[1]$, such that
$g=\si^{-2}\bg$ is Einstein, and parallel standard tractors $\Pa$
with the property that $X_A\Pa^A$ is nowhere vanishing. The mapping
from Einstein scales to parallel tractors is given by 
$\si\mapsto \frac{1}{n}D_A \si$ while the inverse is $\Pa^A \mapsto X_A\Pa^A$.
\end{theorem}

From this theorem it is natural to investigate the interpretation of
parallel tractors in general, i.e.\ without assuming $X_A \Pa^A$ is
non-vanishing.  If $\nd_a \Pa_B=0$ at $x\in M$, then at $x$ we have
$\Pa_B =\frac{1}{n}D_B \si$.  So for parallel $\Pa^A$, it follows from
the definition of $D$ that $X_A \Pa^A$ vanishes on an open set if and
only if $\Pa^A$ vanishes on the same open set. So suppose $\Pa^A\neq
0$ is parallel (and hence non-vanishing).  Then $\si:=X_A \Pa^A$ may
vanish at points or on suitable closed sets, but the 2-jet of $\si$ is
non-vanishing.  In particular, $\si$ is non-vanishing on an open dense
set. From this observation and the Theorem, it follows that 
a manifold admits a parallel tractor if and only if there is a section
$\si$ of $\ce[1]$ such that on some open dense subset $\si^{-2}\bg$ is
Einstein. Since, for a tractor $\Pa$, $\si=X_A\Pa^A$
determines a metric by $g=\si^{-2}\bg$ let us call any point where $\si$
vanishes a point of scale singularity.

Let us say a metric $g$ (or a conformal structure $[g]$) on a manifold
$M$ is {\em almost conformally Einstein} if it is conformally Einstein
on an open dense subset of $M$.  The Poincar\'e metrics of
\cite{FeffGr} are examples of such metrics and it is straightforward
to construct other examples \cite{GLeit}.  This is a natural
generalisation of the Einstein condition since any natural conformally
invariant tensor which vanishes for conformally Einstein metrics also
vanishes, by dint of continuity, for almost conformally Einstein
metrics. The requirement that a metric admit a parallel standard
tractor (which is ostensibly a stronger requirement than being almost
conformally Einstein) is an especially natural condition from two
points of view. Firstly although it is a conformal condition each
parallel tractor determines an actual Einstein metric on an open dense
subset of the manifold. Secondly the equation of parallel transport
for the tractor field gives a canonical extension of the Einstein
equations through any scale singularities. Hence a Riemannian or
pseudo-Riemannian structure (or a conformal structure $[g]$) equipped
with a parallel standard tractor will be termed an {\em almost
  Einstein} structure.

Let us make a few simple observations concerning almost conformally
Einstein structures. Suppose that $(M,[g])$ is an almost conformally
Einstein manifold such that the (maximal) conformally Einstein subset
is connected.  Then it has a parallel tractor on a connected open
dense subset of $M$ and thus (extending this by parallel transport) a
parallel tractor on $M$. For example Fefferman and Graham's Poincar\'e metrics
admit a parallel tractor and so are almost Einstein.
Note that when we have a parallel tractor $\Pa$, then the constant
function $-\frac{n}{2}\Pa^A\Pa_A$
extends the scalar curvature $\J$ across points of scale singularity.

Since $\Pa$ parallel implies $\Pa=\frac{1}{n}D \si$ for some density
of weight 1, an obvious question is whether, at any point $x$, we can
have $j^1_x\si=0$. In fact this can happen since on conformally flat
structures the tractor connection is flat; we may set $\Pa=(0,0,1)$ at
some point $x$ and then extend (locally at least) to a parallel field
by parallel transport. We will see below (see the comment following
Theorem \ref{weaklynot}) that if $\Pa=\frac{1}{n}D \si$ (non-trivial)
is parallel and $j^1_x\si=0$ then the Weyl curvature must vanish at
$x$. On top of this there there are other severe restrictions.
\begin{proposition}\label{j1facts}
If $\Pa=\frac{1}{n}D \si\neq 0$ is parallel
and $j^1_x\si=0$ then $\si\neq f_1f_2 \tilde{\si}$ where $\tilde{\si}$ is a section of $\ce[1]$, and  $f_1$ and
$f_2$ are functions which vanish at $x$. 

If $\Pa=\frac{1}{n}D \si\neq 0$ is parallel then $j^1\si$ can only
vanish at isolated points. In Riemannian signature if
$\Pa=\frac{1}{n}D \si\neq 0$ is parallel and $j^1_x\si=0$ then 
there is a neighbourhood of $x$ such that, in this
neighbourhood, $\si$ is non-vanishing away from  $x$.
\end{proposition}
\noindent{\bf Proof:} For the first statement suppose that
$\Pa=\frac{1}{n}D \si\neq 0$ is parallel and $j^1_x\si=0$. Then
(without loss of generality, by a constant rescaling we may assume)
$\Pa (x)=(0,0,-1)$. Thus from \nn{trconn} we have $\nd_a \nd_b \si
(x)=g_{ab}(x)$ (whereas $\nd_a \nd_b (f_1f_2 \tilde{\si})$ has rank at
most 2 at $x$).

For the second claim suppose  that $\Pa=\frac{1}{n}D \si\neq 0$ is parallel 
 and $j^1\si$ vanishes
along a curve $x(t)$ through $x$. Then if $t^a$ is the tangent field
we have $t^a\nd_a \nd_b \si=0$ along the curve and in particular at
$x$. But at $x$, as observed, we have $\nd_a\nd_b\si=g_{ab}$ and so
$t^a\nd_a \nd_b \si= t^ag_{ab}\neq 0$ and we have a contradiction.

Finally let us pick a metric $g$ from the conformal class and with
this trivialise the density bundles.  If $\Pa=\frac{1}{n}D \si\neq 0$
is parallel and $j^1_x\si$ vanishes then $\nd_a\nd_b \si
(x)=g_{ab}(x)$ and so, in terms of Riemann normal coordinates based at
$x$, the first non-vanishing term in the Taylor series for $\si$
(based at $x$) is $g_{ij}x^ix^j$. 
\quad $\Box$

\section{Obstructions}\label{three}

There is an obvious integrability condition for the existence of a parallel tractor:
$$ 
 \nd_a\Pa^C=0 \Rightarrow [\nd_a,\nd_b]\Pa^C=\Om_{ab}{}^C{}_D\Pa^D=0
$$
  where $\Om_{ab}{}^C{}_D$ is the curvature of
the tractor connection.  So $\Om_{ab}{}^C{}_D\Pa^D=0$ is a necessary
condition for a metric to be conformal to Einstein.
It is an elementary exercise to interpret this in terms of tensors.
An easy calculation establishes that
\begin{equation}\label{tractcurv}
\Om_{ab}{}^C{}_D=
\left(\begin{array}{lll}
0&0&0\\
A^{c}{}_{ab} & C_{ab}{}^c{}_d & 0\\
0 & -A_{dab} & 0
\end{array}\right)
\end{equation} 
where $A_{abc}:=2\nabla_{[b}\Rho_{c]a} $ is the {\em Cotton tensor}
for the metric $g$ giving the splitting of the tractor
bundles. 
Thus, if 
$$\left(\begin{array}{c}\si\\ \mu^d \\\rho \end{array}\right)
$$  
gives $\Pa^D$ in the scale $g$   then 
\begin{equation}\label{Ctilde}
 \quad\quad \Om_{ab}{}^C{}_D\Pa^D=0 \Leftrightarrow  \si A^c{}_{ab}+\mu^d C_{ab}{}^c{}_d=0 .
\end{equation}
(For the implication $\Leftarrow$ note that $\si A^c{}_{ab}+\mu^d
 C_{ab}{}^c{}_d=0$ implies $\mu^d A_{dab}=0$ from the symmetries of
 the Weyl curvature.)  
Where $\si$ is non-vanishing this is
$$ 
\mbox{[C]} \quad \quad A_{dab}+K^c C_{cdab}=0
$$ with $K_c=-\si^{-1}\mu_c$. If we require in addition that $K_c$ is
exact then this is the C-space equation of \cite{Sz,CTP}, the equation
that is satisfied if and only if the Cotton tensor vanishes in some
scale; $\Om_{ab}{}^C{}_D\Pa^D=0$ generalises this. 

We can now differentiate  to obtain more
integrability conditions. $\Om_{ab}{}^C{}_D\Pa^D=0$ implies
$$
0=\Pa^D\nd_e\Om_{ab}{}^C{}_D + \Om_{ab}{}^C{}_D\nd_e \Pa^D
$$
so if $\Pa^D$ is parallel then 
$$
\mbox{[D]} \quad \quad \Pa^D\nd_e\Om_{ab}{}^C{}_D=0~ .
$$
Obviously we could continue along these lines.  More abstractly the
situation is this. If $\Pa$ is parallel then it is obviously fixed by
the conformal holonomy group for the tractor connection.  Thus the
infinitesimal holonomy group fixes $\Pa$ and the Lie algebra of this
is generated by the tractor derivatives of the tractor curvature. In
general calculating even the infinitesimal holonomy is not possible.
We shall see that in a rather general setting it suffices,for our
purposes, to stop at [D]. Note that the pair $(\Om_{bc}{}^D{}_E,\nd_a
\Om_{bc}{}^D{}_E)$ simply recovers the 1-jet of the curvature at each
point, and so the pair is conformally invariant.

Contracting [D] via $g^{ea}$ gives $ \Pa^D\nd^a\Om_{ab}{}^C{}_D=0 $
and via a short calculation we may re-express this by the equation
$$
\mbox{[B]} \quad \quad B_{ab}+(n-4)K^dK^c C_{dabc}=0 ,
$$
provided $X_A\Pa^A$ is non-vanishing.
Here 
$$
B_{ab}:=\nabla^c
A_{acb}+\Rho^{dc}C_{dacb}.
$$
is the {\em Bach tensor} (which is conformally invariant in
dimension 4).  Since the pair $(\Om_{bc}{}^D{}_E,\nd_a
\Om_{bc}{}^D{}_E) $ is conformally invariant it is clear that so is
the system [C], [D] and hence so also the system [C], [B]. In
dimension 4 the latter is exactly the system considered by Kozameh,
Newman and Tod in \cite{CTP}. In any dimension assuming the metric is
sufficiently generic this generalisation of the KNT system is
sufficient to determine whether or not a metric is conformally
Einstein (see theorems 2.2 and 2.3 of \cite{GNur}). 

For the purposes of characterising conformally Einstein metrics or
almost conformally Einstein metrics the system [C], [D] has more
information. Let us say that a (pseudo-)Riemannian manifold is {\em
  weakly generic} if, at each $x\in M$, the only solution $V^d_x\in
T_xM$ to
$$
C_{abcd}V^d_x=0~, \quad \mbox{ at } x\in M~,
$$
is $V^d_x=0$. 
We will say that a (pseudo-)Riemannian manifold is {\em almost 
weakly generic} if the only smooth vector field  $V^d$ solving
$$
C_{abcd}V^d=0~, 
$$
is $V^d=0$.
The following is a  generalisation of Theorem 3.4 of \cite{GNur}. 
\begin{theorem} \label{tracchar}
 A weakly generic conformal manifold is conformally 
 Einstein if and only if there exists a non-vanishing tractor field
 $\Pa^A\in \ce^A$ such that
$$
\begin{array}{cl}
\Pa^E\Om_{bcDE}=0 & \mbox{\rm [$\tilde{\rm C}$]} \\
\Pa^E\nd_a\Om_{bcDE}=0 & \mbox{\rm [$\tilde{\rm D}$]}. 
\end{array}
$$

An almost weakly generic conformal manifold admits a parallel standard
tractor if and only if there exists
a non-vanishing tractor field $\Pa^A\in \ce^A$ such that \mbox{\rm
  [$\tilde{\rm C}$]} and \mbox{\rm [$\tilde{\rm D}$]} hold.
\end{theorem}
\noindent {\bf Proof:} $\Rightarrow :$ is clear for both statements. \\
$\Leftarrow :$ We treat the first statement first. So let us assume
the structure is weakly generic.  Let $(\si,\mu^e,\tau)$ be the
components of $\Pa^E$. Following \cite{GP} we write $Y^E $ and
$Z^E{}_e$ for the injectors (determined by the choice metric giving
the splitting) into the first two tractor slots. Since 
$[\ce^A]_g=\ce[1]\oplus \ce_a[1]\oplus \ce[-1] $ we have  
$$
\hspace*{-10mm} X^A:\ce[-1]\to \ce^A \quad Z^{Aa}:\ce_a[1]\to \ce^A 
\quad Y^A:\ce[1]\to \ce^A
$$
so that we may write
$\Pa^E= Y^E\si + Z^{Ed}\mu_d+X^E\t$.  Suppose that $X_A\Pa^A=\si$
vanishes at some point $x$.  Then from \nn{Ctilde} we have $ \mu^d
C_{abcd}=0$ at $x$ and so, since the conformal class is weakly
generic, $\mu^d(x)=0$. Thus $\Pa^E=\t X^E$, at $x$, and [$\tilde{\rm
D}$] gives $X^E\nd_a\Om_{bcDE}=0$ at $x$. But from \nn{trconn} $\nd_a
X^E=Z^E{}_a$ and from \nn{tractcurv} $X^E\Om_{bcDE}=0$, and so $
Z_D{}^d C_{bcda}-X_{D} A_{abc}= Z^E{}_a\Om_{bcDE}=0$ at $x$. But this
means $C_{bcda}(x)=0$ which contradicts the assumption that the
conformal class is weakly generic. So $X_A\Pa^A$ is non-vanishing.

Differentiating [$\tilde{\rm C}$] and using
[$\tilde{\rm D}$] gives
$$
\Om_{ab}{}^C{}_D\nd_e\Pa^D=0.
$$ Since $g$ weakly generic $\Om_{bc}{}^{D}{}_{E}$ must have rank at
least $n$ as a map $\Om_{bc}{}^D{}_E:\ce^{bc}{}_D\to \ce_E$. On the
other hand [$\tilde{\rm C}$] and the formula for $\Om_{bc}{}^D{}_E$ we
have that $\Pa^E$ and (the linearly independent) $X^E=(0,0,1)$ are
orthogonal to its range. So
$$
\nd_a\Pa^E= \alpha_a \Pa^E +\beta_a X^E ,
$$
for some 1-forms $\alpha_a$ and $\beta_a$. Differentiating again and 
alternating 
leads to $\nd_{[a}\alpha_{b]}=0$. So locally $\alpha =df$ (some function $f$) and 
$\tilde{\Pa}^E:= e^{-f}\Pa^E$ satisfies
$$
\nd_a\tilde{\Pa}^E= \tilde{\beta}_a X^E , \quad \mbox{that is:}
$$
$$
\begin{array}{l} \nabla_a \tilde{\si}-\tilde{\mu}_a =0 \\
                       \nabla_a \tilde{\mu}_b+ \bg_{ab} \tilde{\tau} +\V_{ab}\tilde{\si} =0 ~.
   \end{array} 
$$ 
So, calculating in terms of the metric $g:= \tilde{\si}^{-2}\bg$, we have
$\tilde{\mu}_a =\nabla_a \tilde{\si}=0 $ and $\V_{ab} + \bg_{ab}
\tilde{\tau}/\tilde{\si}=0$. That is the metric $g$ is Einstein and
$\frac{1}{n}D_A \tilde{\si}$ is parallel.  

Now the implication $\Leftarrow$ follows easily for the second
statement. If the structure is almost weakly generic then it follows
easily that it is weakly generic on an open dense subset of $M$. Thus
we obtain a parallel tractor on an open dense subset which then
extends by continuity.  \quad $\Box$

We have the following consequence of the theorem above.
\begin{corollary}\label{weakanysig}
An weakly generic (almost weakly generic) pseudo-Riemannian or
Riemannian metric $g$ on an $n$-manifold is conformally Einstein
(resp.\ admits a parallel standard tractor) if and only if the natural
invariants
$$
\Om_{abK D_1}
 \cdots \Om_{cdLD_s} \nd_{e}\Omega_{fgPD_{s+1}} \cdots \nd_h \Omega_{k\ell QD_{n+2}} ,
$$ for $s=0,1,\cdots,n+1$, all vanish identically. Here the sequentially labelled
indices $D_1,\cdots , D_{n+2}$ are completely skewed over.
\end{corollary}
\noindent{\bf Proof:} The Theorem can clearly be rephrased to state
that $g$ is conformally Einstein (or, simply admits a parallel tractor, in
the almost weakly generic case) if and only if the map
\begin{equation}\label{Ommap}
(\Om_{bcDE},\nd_a\Om_{bcDE}):\ce^{bcD}\oplus\ce^{abcD} \to \ce_E
\end{equation}
given by 
$$
(V^{bcD},W^{abcD})\mapsto V^{bcD} \Om_{bcDE}+W^{abcD}\nd_a\Om_{bcDE}
$$ fails to have maximal rank at every point of $M$. But by elementary
linear algebra this happens if and only if the induced alternating
multi-linear map to $\Lambda^{n+2}(\ce^E)$ vanishes. This is
equivalent to the claim in the Corollary, since for any metric the
tractor curvature satisfies $\Om_{bcDE}X^E=0$.  \quad $\Box$ \\
The power of Corollary \ref{weakanysig} is that it gives sharp
curvature obstructions for any signature and the idea clearly
generalises in an obvious way to related problems on other
structures. 

We should point out that in the case of Riemannian signature there is
a much simpler sharp obstruction to conformally Einstein metrics in
the weakly generic setting.
Let us  write $L^a_b:= C^{acde}C_{bcde}$
and $\tilde{L}^a_b$ for the tensor field which is the
pointwise adjugate of $L^a_b$.  $\tilde{L}^a_b$ is given by a formula
which is a partial contraction polynomial (and homogeneous of degree
$2n-2$) in the Weyl curvature and for any structure we have
$$
\tilde{L}^a_b L^b_c =||L||\delta^a_c ,
$$ where $||L||$ denotes the determinant of $L^a_b$.
Now  define 
$$
D^{acde}:=- \tilde{L}^a_b C^{bcde} 
$$ Then $D^{acde}$ is a natural conformal invariant and if [C] has a
solution for $K^e$ (which is necessary for the metric to be
conformally Einstein) then by uniqueness the solution is
$||L||^{-1}D^{edab}A_{dab}$. If $K_e$ is exact then it is the
$\Upsilon_e$ in \nn{confein} and so we have the following result
\cite{GNur}.
\begin{theorem}\label{nob}
The natural invariant 
$$
G_{ab}:=\mbox{Trace-free}\big[||L||^2 \Rho_{ab}-||L||\nd_a(
D_{bcde}A^{cde})+\quad\quad\quad\quad\quad\quad\quad\quad\quad\quad\quad\quad $$
$$\quad\quad\quad\quad\quad\quad\quad\quad\quad\quad\quad\quad(\nd_a||L||)(
D_{bcde}A^{cde}) +D_{aijk}A^{ijk} D_{bcde}A^{cde} \big].
$$ 
 is a conformal invariant of weight
  $-8n$.  A manifold with a weakly generic Riemannian metric $g$ is
  conformally Einstein if and only if $G_{ab}$ vanishes. The same is
  true on pseudo-Riemannian manifolds where the conformal invariant
  $||L||$ is non-vanishing.
\end{theorem}

We observe here that weakly generic almost conformally
Einstein metrics do not admit parallel tractors with scale singularities.
This observation motivates the consideration above of
 almost weakly generic metrics.
\begin{theorem}\label{weaklynot}
If a weakly generic conformal structure admits a non-zero parallel
standard tractor field $\Pa^A$ then it is conformally Einstein and
$X_A\Pa^A$ is an Einstein scale.
\end{theorem}
\noindent{\bf Proof:} If $\Pa^A$ is parallel then the system
[$\tilde{\rm C}$] and [$\tilde{\rm D}$] holds. Early in the proof of
Theorem \ref{tracchar} it is established that if these hold and the
structure is weakly generic then $X_A\Pa^A$ is non-vanishing.
\quad $\Box$ \\ 
Finally a related observation concerning the analysis
of points of scale singularities. 
 If $\Pa_A$ is parallel then, recall,
$\Pa_A=\frac{1}{n}D_A \si$ for some density of $\si$ of weight 1. If
at some point $x$ we have $j^1_x\si$ then at $x$ we have $\Pa=\tau
X^A$ and so from [$\tilde{\rm D}$] and arguing once again as in the
proof of Theorem \ref{tracchar} we have that if $C_{abcd}(x)=0$.

\section{The Fefferman-Graham ambient obstruction tensor}

Recall that the splitting of the tractor bundle given by a metric is
determined by the formula $( \si, \nd_a \si, - \frac{1}{n}(\Delta +
\J) \si )$ for the operator $\frac{1}{n} D$ on densities of weight 1.
From this it is easily verified that, under a conformal transformation
$g\mapsto \wh{g}=e^{2\om} g$, the injector $Z_A{}^a$ transforms to
$Z_{A}{}^a+\Up^aX_A$. With the convention that sequentially
labelled subscript indices, e.g.\ $A_0$, $A_1$, and $A_2$, are
implicitly skewed over, it follows immediately that
$$
X_{A_0}
Z_{A_1}{}^bZ_{A_2}{}^c\Omega_{bc}{}^D{}_E,
$$
is conformally invariant.  The operator $D$ extends to an operator
(the so-called tractor-D operator) on densities and tractors of
general weight $w$ by the formula
$$ 
D^A V:=(n+2w-2)w Y^A V+ (n+2w-2)Z^{Aa}\nabla_a V -X^A\Box V,
$$
where $\Box V:= \Delta V+ w \J V$ and in these formulae $\nd$ is the
coupled Levi-Civita-tractor connection. Using the conformal
transformation formula for $Z_A{}^a$ and the corresponding formula for
$Y^A$ it is easily verified directly that this is conformally
invariant. 
Thus, by construction the W-tractor field 
$$
W_{A_1A_2}{}^D{}_E := \frac{3}{n-2}D^{A_0}X_{A_0}
Z_{A_1}{}^bZ_{A_2}{}^c\Omega_{bc}{}^D{}_E,
$$ is conformally invariant. It is readily verified that it has Weyl
curvature type symmetries:
$$
W_{ABCD}=W_{[AB][CD]}, \quad  W_{[ABC]D}=0 \quad \mbox{Trace-free},
$$
etcetera.

Recall that if $\Pa^E$ is a parallel tractor then
$\Omega_{bc}{}^D{}_E\Pa^E=0$ and from this and the definition of the W-tractor 
we have 
\begin{equation}\label{WPa}
W_{A_1A_2}{}^D{}_E \Pa^E=0,
\end{equation}
since (viewed as an order 0 operator) $\Pa^E$ commutes with the
tractor connection and, indeed, the tractor-D operator. So the existence of a
tractor satisfying \nn{WPa} is a necessary condition for a metric to
admit a parallel tractor and, in particular, for a metric to be
conformally Einstein.  In fact in dimensions other than 4, the
equation \nn{WPa} is in fact exactly an alternative expression for the
system [C], [B] (see \cite{GNur}).

The situation is more interesting in dimension 4. Expanding $W_{ABCE}$
out we obtain,
$$
W_{ABCE}=\left(\begin{array}{ccccc}
&&0&& \\
&0&&0& \\
(n-4)C_{abce}& & 0\hspace*{5mm} 0 && 0\\
& (n-4)A_{bce} &&0&\\
&& B_{eb}&&
\end{array}\right)
$$
where we have indicated components with respect to the composition
series for tractors (of weight $-2$) with Weyl tensor type symmetries
(see \cite{GoPetob} for details). So in dimension 4 the only possibly
non-vanishing component of the W-tractor is the Bach tensor. This is
extracted by the tractor for a scale $\si$.  For calculations it is
more convenient to express the expansion in terms of the
injectors/projectors $X$, $Y$ and $Z$. In dimension $n$ we have that
$W_{ABCE}$ is given by
\begin{equation}\label{Wform}
\begin{array}{l}
(n-4)\left( Z_A{}^aZ_B{}^bZ_C{}^cZ_E{}^e C_{abce}
-2 Z_A{}^aZ_B{}^bX_{[C}Z_{E]}{}^e A_{eab}\right. \\ 
\left.-2 X_{[A}Z_{B]}{}^b Z_C{}^cZ_E{}^e A_{bce} \right)
+ 4 X_{[A}Z_{B]}{}^b X_{[C} Z_{E]}{}^e B_{eb}. 
\end{array}
\end{equation}
Let $\si$ be a choice of conformal scale and $\Pa_E:=\frac{1}{n}D_E
\si$.  Splitting the standard tractor bundle via the metric
$g=\si^{-2}\bg$ given by the scale $\si$ (and noting that for the
connection given by this scale we have $\nd\si=0$) we have $
\Pa_E=(\si, 0, -\frac{1}{n}\J \si), $ or in other terms $\Pa_E= \si
Y_E -\frac{1}{n}\si\J X_E$. Contracting this into the last display and
setting $n=4$ we obtain
$$ 
W_{ABCE} \Pa^E= -2 \si X_{[A}Z_{B]}{}^b Z_{C}{}^e B_{eb}. 
$$ (For the contraction note that from the formula for the tractor
metric we have that $X^E$, $Y^E$ are null, $X^EY_E=1$, and both are
orthogonal to $Z_E{}^e$.) Now if $\si$ is an Einstein scale (i.e.\
$g=\si^{-2}\bg$ an Einstein metric) then this vanishes by \nn{WPa}. Of
course there are easier ways to show that the Bach tensor is an
obstruction to conformally Einstein metrics. The main point is that
this generalises.

Note that the W-tractor has conformal weight $-2$. In dimension 6 this
is $1-n/2$ and so exactly the weight on which the tractor coupled
conformal wave operator $\Box$ (as defined earlier) acts invariantly.  In
dimension 6 a straightforward direct calculation shows that
$$
\fb W_{ABCE}= \left(\begin{array}{ccccc}
&&0&& \\
&0&&0&  \\
0& & 0 \hspace*{5mm} 0 && 0\\
& 0 &&0&\\
&& \cB^{(6)}_{eb}&&
\end{array}\right)
$$
where $\fb$ is the modification of $\Box$ given by 
$$
\fb=\Box + \frac{1}{4}W\sharp\sharp ,
$$ and $\hash$ indicates the natural action of sections of ${\rm
End}(\ce^A)$ on tractors. (Note that $W$ is a section in
$\otimes^2 ({\rm End}(\ce^A))$ and that if an element of ${\rm
End}(\ce^A)$ is skew relative to the tractor metric $h$ then its action by 
$\hash$  obviously preserves
the symmetry type of tractor bundles.)
  Evidently $\cB^{(6)}$ is some
conformal invariant and from the composition series concerned we know
this is a trace-free symmetric density-valued 2-tensor of conformal
weight $-4$. By essentially the same argument as for the Bach tensor
above we see this an obstruction to conformally Einstein metrics. On
the one hand for any scale $\si$ contracting $\Pa^E:=\frac{1}{n}D^E
\si$ into $\fb W_{ABCE}$ extracts (the non-vanishing scale $\si$
times) $\cB^{(6)}_{eb} $. On the other hand if $\si$ is an Einstein
scale then $\Pa$ is parallel and so commutes with $\fb$. Thus from
\nn{WPa} $\cB^{(6)}_{eb} $ necessarily vanishes for conformally Einstein metrics.

More generally we have the following (an adaption of part of Theorems
4.1 and 4.2 from \cite{GoPetob}).
\begin{theorem}\label{obstrn}
 Let $M$ be a conformal manifold of
dimension $n$ even.  There is a (conformally invariant) operator
$$
\fb_{n/2-2}: \otimes^2(\Lambda^2\ce^A)[-2]\to \otimes^2(\Lambda^2\ce^A)[2-n],
$$ which may be expressed by a formula polynomial in
 $X,\nd_A:=Z_A{}^a\nd_a,W,h$, and $h^{-1}$, such that
$$
\fb_{n/2-2} W_{A_1A_2B_1B_2} =
 K(n) X_{A_1}Z_{A_2}{}^aX_{B_1}Z_{B_2}{}^b \cB^{(n)}_{ab},
$$
 Here
 $K(n)$ is a (known) non-zero constant depending on $n$. In the polynomial formula for the left-hand-side the free indices 
always appear on a $W$. 
\end{theorem}

From this we have immediately the generalisation of the result for the
Bach tensor and its dimension 6 analogue. Once again from an easy
analysis of the composition series for the tractor bundle (with Weyl
tensor symmetries) that $\fb_{n/2-2} W_{A_1A_2B_1B_2}$ takes values in
we can conclude that $ \cB^{(n)}_{ab}$ is a trace-free symmetric
density-valued 2-tensor of conformal weight $2-n$. By construction it
is clearly natural (i.e.\ can be given by a formula polynomial in the
conformal metric and its inverse and the Levi-Civita covariant
derivatives of the Riemannian curvature).
\begin{corollary} 
The natural tensor $\cB^{(n)}_{ab}$ is conformally invariant and
vanishes for almost conformally Einstein metrics.
\end{corollary}
\noindent {\bf Proof:} The conformal invariance is immediate from the
conformal invariance of the operator $\fb_{n/2-2}$, the tractor field
$W_{A_1A_2B_1B_2}$ and the result that $X_{A_1}Z_{A_2}{}^a:\ce_a\to
\ce_{A_1A_2}$ is injective and conformally invariant.

For a scale $\si$ let $\Pa_A:=\frac{1}{n}D_A\si$. Then, in terms of the 
metric $g=\si^{-2}\bg$, we have 
$$ \Pa^A = \si Y^A -\frac{1}{n}\J \si X^A
$$
and so from the theorem (and the tractor metric) we have
$$
\begin{aligned}
\si^2 \cB^{(n)}_{ab}& =
4(K(n))^{-1}Z^{A_2}{}_aZ^{B_2}{}_b\Pa^{A_1}\Pa^{B_1}\fb_{n/2-2}
W_{A_1A_2B_1B_2}.
\end{aligned}
$$
But if $\si$ is an Einstein scale then $\Pa^A$ is parallel and so 
commutes with all terms  in the polynomial expression for $\fb_{n/2-2}
W_{A_1A_2B_1B_2} $ and the free indices $A_1A_2B_1B_2$ always appear
on a $W$. Thus, 
$$
\hspace*{-.5cm}\Pa^{B_1} \fb_{n/2-2} W_{A_1A_2B_1B_2}=0 \quad
\mbox{since} \quad \Pa^{A}W_{ABCD}=0
$$
and so from the previous display it follows that, in the scale $\si$,
$$
\cB^{(n)}_{ab}=0.
$$ Since $\cB^{(n)}_{ab}$ is conformally invariant it must vanish for
any metric which is conformally Einstein and then by continuity for
any metric which is almost conformally Einstein.
\quad $\Box$

The proof of Theorem \ref{obstrn} in \cite{GoPetob} uses the conformal
ambient metric construction of Fefferman-Graham and it follows
easily from this proof that $\cB^{(n)}$ is the obstruction to the
ambient metric construction found by Fefferman-Graham.


\begin{thebibliography}{99}


\bibitem{BEG} T.N. Bailey, M.G. Eastwood, and A.R. Gover,
{\em Thomas's structure bundle for conformal, projective and related
structures}, Rocky\ Mountain\ J.\ Math., {\bf 24} (1994) 1191--1217.


\bibitem{B1} H.W.\ Brinkman, {\em Riemann spaces conformal to Einstein
  spaces} Math.\ Ann.\ {\bf 91} (1924) 269--278.

\bibitem{B2} H.W.\ Brinkman, {\em Einstein spaces which are mapped
    conformally on each other} Math.\ Ann.\ {\bf 94} (1925) 119--145.

  
\bibitem{luminy} A.\ \v Cap and A.R. Gover, {\em Tractor bundles
    for irreducible parabolic geometries}, Global analysis and
  harmonic analysis (Marseille-Luminy, 1999) 129-154, S\'{e}min.\ 
  Congr., {\bf 4}, Soc.\ Math.\ France, Paris 2000.  Preprint ESI 865,
  available for viewing on the internet at {\tt http://www.esi.ac.at}

\bibitem{Cap-Gover} A.\ \v Cap and A.R. Gover, {\em Tractor calculi for
parabolic geometries}, Trans.\ Amer.\ Math.\ Soc., {\bf 354} (2002)
1511-1548. 


\bibitem{CapGoFeff} A. \v Cap and A.R. Gover, in progress. 

\bibitem{Cartan} E. Cartan, \textit{Les espaces \`a connexion
conforme}, Ann. Soc. Pol. Math., \textbf{2} (1923), 171--202.



\bibitem{FeffGr} C. Fefferman and C.R. Graham,
{\em Conformal invariants}, in {\em Elie Cartan et les math\-\'ematiques
d'aujourd'hui}, Ast\'erisque, hors s\'{e}rie (Soci\'{e}t\'{e}
Math\'{e}matique de France, Paris, 1985) 95--116.


\bibitem{Gau} P.\ Gauduchon, {\em Connexion canonique et structures de
    Weyl en geometrie conforme}, Report: CNRS UA766 (1990).


\bibitem{Goadv} A.R.\ Gover, {\em Invariant theory and calculus for
    conformal geometries}, Adv.\ Math., {\bf 163} (2001) 206--257.

\bibitem{GLeit} A.R.\ Gover and F.\ Leitner, in progress.

\bibitem{GNur} A.R.\ Gover and P.\ Nurowski,
{\em Obstructions to conformally Einstein metrics in $n$ dimensions}.
  Preprint, math.DG/0405304, http://arXiv.org.

\bibitem{GP} A.R.\ Gover and L.J.\ Peterson, 
{\em Conformally invariant powers of the Laplacian, 
Q-curvature, 
and tractor calculus},
Comm.\ Math.\ Phys., {\bf 235} (2003) 339--378.

\bibitem{GoPetob} A.R.\ Gover and L.J.\ Peterson, {\em The ambient
obstruction tensor and the conformal deformation complex}, 
 Preprint, math.DG/0408229, http://arXiv.org.

\bibitem{GrH} C.R.\ Graham and K.\ Hirachi, {\em The ambient
    obstruction tensor and $Q$-curvature},
arXive preprint: math.DG/0405068,    http://www.arxiv.org

\bibitem{CTP} C. Kozameh, E. T. Newman, K P Tod  {\em Conformal
  Einstein Spaces} GRG {\bf 17}, (1985) 343-352


\bibitem{Listing} M. Listing, {\em Conformal Einstein spaces in
 $N$-dimensions}, Annals of Global Analysis and Geometry, {\bf 20},
 (2001) 183--197.






\bibitem{ot} R. Penrose and W. Rindler,  Spinors and Space-time, vol 1,
Cambridge Univ.\ Press, 1984.


\bibitem{S} J.A. Schouten, Ricci-Calculus, 2nd.\ ed.\ (1954) Springer-Verlag, Berlin.

\bibitem{Sz} P. Szekeres, {\em Spaces conformal to a class of spaces in general relativity},
Proc.\ R.\ Soc.\ London, Ser.\ A, {\bf 274} (1963) 206--212.

\bibitem{Thomas} T.Y.\ Thomas, {\em On conformal geometry},
  Proc.\ Natl.\ Acad.\ Sci.\ USA,
  {\bf 12} (1926) 352--359.


\end{thebibliography}
\end{document}